\documentclass[10pt,english,a4paper]{article}
\usepackage[T1]{fontenc}
\usepackage[utf8]{inputenc}
\usepackage{lmodern}
\usepackage[english]{babel}
\usepackage{amsmath,amsfonts,amssymb,amsthm}
\usepackage{mathrsfs}
\usepackage{pifont}
\usepackage{bm}
\usepackage[shortlabels]{enumitem}
\usepackage{hyperref}
\usepackage{tikz}

\usepackage[nottoc, notlof, notlot]{tocbibind}
\usepackage{tocloft}
\renewcommand{\cfttoctitlefont}{\large\bfseries}
\renewcommand\cftaftertoctitle{\par\noindent\hrulefill\par\vskip 0em} 
\renewcommand{\cftsecfont}{\normalfont}
\renewcommand{\cftsubsecfont}{\normalfont}

\usepackage{geometry}
\usepackage{setspace}
\makeatletter
\g@addto@macro \normalsize {
 \setlength\abovedisplayskip{5pt plus 2pt minus 2pt}
 \setlength\belowdisplayskip{5pt plus 2pt minus 2pt}
}
\makeatother

\DeclareMathOperator{\Hom}{Hom}
\DeclareMathOperator{\Stab}{Stab}
\DeclareMathOperator{\Id}{Id}
\DeclareMathOperator{\car}{char}

\newcommand{\N}{{\mathbb{N}}}
\newcommand{\Z}{{\mathbb{Z}}}
\newcommand{\F}{{\mathcal{F}}}
\newcommand{\bF}{{\mathbb{F}}}
\newcommand{\cB}{{\mathcal{B}}}
\newcommand{\cI}{{\mathcal{I}}}
\newcommand{\cJ}{{\mathcal{J}}}
\newcommand{\cP}{{\mathcal{P}}}
\newcommand{\cT}{{\mathcal{T}}}
\newcommand{\cX}{{\mathcal{X}}}
\newcommand{\triv}{{\mathbf{1}}}

\usepackage{thmtools}
\declaretheoremstyle[
  spaceabove=\topsep, spacebelow=\topsep,
  headfont=\bfseries\scshape,
  headpunct={:},
  notefont=\mdseries, notebraces={(}{)},
  bodyfont=\normalfont,
  postheadspace=1em,
  qed=\qedsymbol
]{mythmstyle}
\declaretheorem[style=mythmstyle,numbered=no,name=Proof]{dem}

\newtheoremstyle{style1}{8pt}{8pt}{\itshape}{}{\bfseries}{:}{.5em}{}
\newtheoremstyle{style2}{8pt}{8pt}{}{}{\bfseries}{:}{.5em}{}

\theoremstyle{style1}
\newtheorem{thm}{Theorem}[section]
\newtheorem{lem}[thm]{Lemma}
\newtheorem{prop}[thm]{Proposition}
\newtheorem{cor}[thm]{Corollary}
\newtheorem{Rq}[thm]{Remark}
\newtheorem*{thintro}{Theorem}

\theoremstyle{style2}
\newtheorem{defn}[thm]{Definition}
\newtheorem*{Not}{Notation}

\definecolor{r1}{RGB}{175,25,25}
\definecolor{g1}{RGB}{130,130,130}
\hypersetup{colorlinks = true,allcolors=r1!80!black}

\begin{document}

\noindent\makebox[\textwidth][c]{\Large Prime ideals of the Burnside ring of a saturated fusion system}
\vspace{0.2cm}
\begin{center} Nicolas Lemoine \end{center}
\vspace{0.3cm}
\begin{abstract}
In a 1969 article, A. Dress described the prime ideals of the Burnside ring of a finite group $G$ and the inclusion relations between them. One may ask whether similar results exist for the Burnside ring of a fusion system $\F$ on a finite $p$-group $S$. The present paper answers positively, providing a typology for the prime ideals of $A(\F)$ and for those of $A(\F)_{(p)}$ as a consequence.
\end{abstract}
\vspace{0.2cm}
\begin{center} Keywords: \textit{saturated fusion systems, Burnside ring, prime ideals.} \end{center}
\vspace{0.2cm}
{\hypersetup{colorlinks = true,linkcolor=black}
\tableofcontents}
\noindent\hrulefill
\vspace{0.6cm}

\section*{Introduction}
\addcontentsline{toc}{section}{Introduction}

For any finite group $G$, the set $A_+(G)$ of isomorphism classes of finite $G$-sets possesses a semiring structure. Taking the Grothendieck group of $A_+(G)$ we obtain a commutative ring called the Burnside ring of $G$, which we denote by $A(G)$. As a $\Z$-module, $A(G)$ has a natural basis given by classes of transitive $G$-sets $[G/H]$, one per class of $G$-conjugation of subgroups $H$ in $G$.

In this paper we consider a finite $p$-group $S$, where $p$ is a fixed prime number, and we are interested in a particular subring of $A(S)$. Let $\F$ be a saturated fusion system over $S$, i.e. a category whose objects are the subgroups of $S$ and whose morphisms are injective group morphisms between those subgroups, such that certain axioms are satisfied (see Section I.2 in \cite{AKO}). An $S$-set $X$ is said to be $\F$-stable if for every subgroup $P$ of $S$ and every morphism $\varphi : P \to S$ in $\F$, the $P$-set induced by restriction of the $S$-action on $X$ along $\varphi$ is $P$-isomorphic to the $P$-set obtained by restriction of the $S$-action on $X$ along the inclusion $P \hookrightarrow S$. Isomorphism classes of $\F$-stable $S$-sets form a subsemiring $A_+(\F)$ of $A_+(S)$, and its Grothendieck group $A(\F)$ is a subring of $A(S)$ that we call the Burnside ring of $\F$. It was proved by S. P. Reeh in his article \cite{Reeh1} that $A_+(\F)$ possesses a $\Z$-basis formed with classes of $\F$-stable $S$-sets, each basis element $\alpha_{[P]_\F}$ being associated to (and constructed from) a unique $\F$-isomorphism class $[P]_\F$ of subgroups of $S$.

In his article \cite{Dress}, published in 1969, A. Dress studied the prime ideals of $A(G)$, giving an explicit decription of those ideals and the inclusion relations between them. The aim of our paper is to exhibit similar results for the prime ideals of $A(\F)$, when $\F$ is a saturated fusion system. For this, we adapt the proof of A. Dress. Our goal is achieved with Theorem \ref{typology} and Proposition \ref{inclusions}, which constitute a generalization of the result of Dress when applied to $p$-groups. The core of these results can be stated as follows, where $\Phi_P : A(S)\to \Z$ denotes the fixed-point homomorphism associated to the subgroup $P$ of $S$ (see Definition \ref{def trace}):

\begin{thintro}
All prime ideals of $A(\F)$ are listed below.
\begin{enumerate}[$\bullet$]
\item \textup{\textbf{"Type $p$" ideal:}} $$\cI_{S,p} = \{ X \in A(\F) \ \vert \ \Phi_S(X) \equiv 0 \mod p \} \ .$$
\item \textup{\textbf{"Type $0$" ideals:}} there is one "type $0$" prime ideal per $\F$-isomorphism class $[P]_\F$, which is $$ \cI_{P,0} = \{ X \in A(\F) \ \vert \ \Phi_P(X) =0 \} \ . $$
\item \textup{\textbf{"Type $q$" ideals, for $q \neq p$ prime:}} there is one "type $q$" prime ideal per class $[P]_\F$, namely $$ \cI_{P,q} = \{ X \in A(\F) \ \vert \ \Phi_P(X) \equiv 0 \mod q \} \ . $$
\end{enumerate}
Moreover, "type $0$" ideals are minimal prime ideals, and "type $q$" ideals for $q$ prime are maximal ideals.
\end{thintro}

We begin the paper by recalling definitions and useful facts about Burnside rings and fusion systems in Section \ref{sec prereq}. Then we give an exhaustive description of the set of prime ideals of $A(\F)$ in Section \ref{sec descri}, that we complete in Section \ref{sec inclusions} with  the inclusion relations between these ideals. The description of prime ideals of the $p$-localized Burnside ring $A(\F)_{(p)}$ comes in Section \ref{sec local}, as a corollary of the study in $A(\F)$.

\subsection*{Acknowledgments}

I owe the idea of adapting this result of A. Dress to Radu Stancu, whom I also thank for his various advice. I am grateful to Rémi Molinier as well, for his helpful suggestions, listening and support. Finally, I thank Jean-Baptiste Meilhan for reviewing this paper before its submission.

\section{Prerequisites}
\label{sec prereq}

The subgroup relationship between groups will be classically denoted by $\leq$. When $G$ is any group, the preorder $\lesssim_G$ on subgroups of $G$ denotes $G$-subconjugation, meaning that $$ P \lesssim_G Q \iff \textup{there exists $ g \in G \ $ such that } \ {}^gP \leq Q \ , \quad \textup{where }  \ {}^gP := gPg^{-1}. $$
When $P \lesssim_G Q$ and $Q \lesssim_G P$, the subgroups $P$ and $Q$ are conjugate in $G$, and we denote this by $P \sim_G Q$. We denote by $Cl(G)$ the set of $G$-conjugacy classes of subgroups of $G$. The class of a subgroup $P \leq G$ in $Cl(G)$ will be denoted by $[P]_G$. The $G$-subconjugation induces a partial order on $Cl(G)$, that we denote by $\lesssim$.

\subsection{The Burnside ring of a finite group}

Let $G$ be a finite group. Let $A_+(G)$ denote the set of isomorphism classes of finite $G$-sets, equipped with the additive law $+$ induced by the disjoint union of $G$-sets, and with the multiplicative law $\times$ induced by the cartesian product of $G$-sets (with diagonal action). This forms a semiring. 

The Burnside ring $A(G)$ of $G$ is defined to be the Grothendieck group of $A_+(G)$. It is commutative and its unit is the class $[G/G]$ of the one-element $G$-set. It also has a natural $\Z$-basis given by classes of transitive $G$-sets $[G/H]$, one per class of $G$-conjugation of subgroups $H \leq G$. The following lemma gives a formula for the product of two transitive $G$-sets. It is elementary and well-known (for example see Formula (3.1) in \cite{Reeh1}), even if its short proof is rarely detailed in papers.
\begin{lem} \label{product formula}
Let $P$ and $Q$ be subgroups of $G$, and let $P \backslash G /Q$ denote the set of double cosets $\overline{g}:=P g Q$ for $g \in G$. We have: $$ [G/P] \times [G/Q] = \sum\limits_{\overline{g} \in P\backslash G/Q} [G/P\cap {}^gQ] $$
\end{lem}

\begin{dem}
Each couple $(g_1P,g_2Q) \in (G/P) \times (G/Q)$ is in the same orbit as $(P,g_1^{-1}g_2Q)$ for the diagonal action of $G$. Moreover, $(P,gQ)$ and $(P,g'Q)$ are in the same orbit if and only if there exists $\widetilde{g} \in G$ such that $\widetilde{g}P=P$ and $g'^{-1}\widetilde{g}g \in Q$, i.e. if and only if $g \in Pg'Q$. To conclude, when $g \in G$ is fixed we can compute the stabilizer subgroup of $(P,gQ)$: \begin{align*}
\Stab_G\big((P,gQ)\big) &= \{ g' \in G \ \vert \ g' \cdot (P,gQ) = (P,gQ) \} \\
&= \{ g' \in G \  \vert \ g' \in P \ \textup{ and } \ g^{-1}g'g \in Q \} \\
&= P \cap {}^gQ
\end{align*}
\end{dem}

In particular, each transitive component of the product $[G/P] \times [G/Q]$ is associated to a subgroup $R$ such that $R \lesssim_G P$ and $R \lesssim_G Q$.

Now we recall the description of the Burnside ring $A(G)$ in terms of fixed points of $G$-sets.

\begin{defn}\label{def trace}
For each subgroup $P \leq G$, we define a $\Z$-linear form $\Phi_P : A(G) \to \Z$ by setting $$ \Phi_P \big([X] \big) := \big|X^P \big| $$ for each isomorphism class of $G$-sets $[X]$, where $X$ is any representative of the class, and where $X^P$ denotes the set of fixed points for the $P$-action on $X$. It is well defined on $A_+(G)$ and extends by linearity to $A(G)$. Since it also preserves products and unit, we obtain a ring homomorphism $\Phi_P : A(G) \to \Z$. Moreover, $\Phi_P$ only depends on the $G$-conjugacy class of the subgroup $P$.
\end{defn}

The value of $\Phi_P$ on transitive $G$-sets is given by the following classical formula (see Section 3 in \cite{Reeh1}). Once again, the proof is simple though often omitted in papers.
\begin{lem} \label{formule phi}
Let $P$ and $Q$ be subgroups of $G$. Then we have
$$ \Phi_P \big([G/Q] \big) = \frac{|N_G(P,Q)|}{|Q|} $$
where $N_G(P,Q)$ denotes the transporter from $P$ to $Q$ in $G$, i.e. the set of all $g \in G$ such that ${}^gP \leq Q$.
\end{lem}

\begin{dem}
A coset $gQ \in G/Q$ belongs to $(G/Q)^P$ if and only if $Pg \subseteq gQ$, i.e. if and only if $g^{-1}Pg \leq Q$, which is equivalent to $g^{-1} \in N_G(P,Q)$. Thus we have a surjective map $N_G(P,Q) \to (G/Q)^P$ given by $g \mapsto g^{-1}Q$, and the preimage of each coset $gQ$ under this map has a cardinality of $|Q|$ (corresponding to the choice of a representative $g$ for $gQ$). We then deduce that $$ \big|N_G(P,Q)\big| = \big|(G/Q)^P\big| \cdot |Q| \ .$$ 
\end{dem}

In particular, $\Phi_P \big([G/Q] \big) \neq 0$ if and only if $P \lesssim_G Q$.

\begin{defn}
We define the mark homomorphism, also called ``ghost map'' in the literature, $$ \Phi : A(G) \longrightarrow \prod\limits_{[P] \in Cl(G)} \Z $$ to be the product of the morphisms $\Phi_P$ for $[P]_G \in Cl(G)$, i.e. for each $a \in A(G)$ we have: $$\Phi(a) = \big(\Phi_P (a) \big)_{[P]_G \in Cl(G)}$$
\end{defn}

This defines a ring homomorphism, and W. Burnside proved that it is injective with finite cokernel (see Section 181 in \cite{Burn}). Hence it provides a practical description of $A(G)$, as a subring of $\Z^d$ where $d=|Cl(G)|$.

\subsection{Fusion systems}

Let $S$ denote a $p$-group, where $p$ is a fixed prime number.

\begin{defn} \label{def fusion system}
A fusion system $\F$ over $S$ is a (small) category whose objects are the subgroups of $S$ and whose morphism sets $\F(P,Q)$, for $P$ and $Q$ subgroups of $S$, satisfy the following axioms:
\begin{enumerate}[$(i)$]
\item Every morphism in $\F(P,Q)$ is an injective group homomorphism.
\item The restriction at $P$ of every inner morphism $c_s$, with $s \in S$ such that $^sP \leq Q$, belongs to $\F(P,Q)$.
\item If $\varphi : P \to Q$ is a morphism of $\F$, then both $\varphi : P \to \varphi(P)$ and $\varphi^{-1} : \varphi(P) \to P$ are morphisms of $\F$.
\end{enumerate}
\end{defn}

The crucial example of a fusion system over $S$ is the one induced by a supergroup $G$ containing $S$, in which the morphisms are the restrictions of inner morphisms of $G$ to subgroups of $S$ in $G$. This fusion system is denoted by $\F_S(G)$. In particular, $\F_S(S)$ is called the inner fusion system of $S$, and it is minimal among fusion systems over $S$ because every fusion system over $S$ contains the morphisms in $\F_S(S)$, according to Axiom $(ii)$ in Definition \ref{def fusion system}.

In our paper we will consider only saturated fusion systems, which are fundamental instances of fusion systems that generalize $\F_S(G)$ in the case where $S$ is a $p$-Sylow subgroup of $G$ (see Theorem 2.3 in Part I of \cite{AKO}). The definition is obtained by adding certain axioms to Definition \ref{def fusion system}. We omit it here, since we will only use saturation through Proposition \ref{base alpha}, and we refer the reader to Section I.2 in \cite{AKO} for a precise definition and details.

For a fusion system $\F$ over a finite $p$-group $S$, we will denote the set of morphisms from $P$ to $Q$ in $\F$ by $\F(P,Q)$. The preorder $\lesssim_\F$ on subgroups of $S$ denotes $\F$-subconjugation, meaning that $$ P \lesssim_\F Q \iff \exists \varphi \in \F(P,Q) \iff \exists \varphi \in \F(P,S), \quad \varphi(P) \leq Q \ .$$
When $P \lesssim_\F Q$ and $Q \lesssim_\F P$, the subgroups $P$ and $Q$ are $\F$-isomorphic (or $\F$-conjugate), and we denote this by $P \sim_\F Q$. We denote the set of $\F$-isomorphism classes of subgroups of $S$ by $Cl(\F)$. The class of a subgroup $P \leq S$ in $Cl(\F)$ will be denoted by $[P]_\F$. The $\F$-subconjugation between subgroups of $S$ induces a partial order on $Cl(\F)$, that we denote by $\lesssim$.

\subsection{The Burnside ring of a fusion system}

Now we will consider a finite $p$-group $S$ and fusion system $\F$ over $S$. If $T$ is another finite group and $X$ is an $S$-set, then any homomorphism of groups $\varphi : T \to S$ induces a $T$-action on $X$ via: $$ \forall x \in X, \ \forall t \in T, \qquad t \cdot x := \varphi(t) \cdot x $$
This gives a well-defined semiring homomorphism $A_+(S) \to A_+(T)$, which extends to a ring homomorphism $r_{\varphi} : A(S) \to A(T)$. In particular when $P \leq S$ and $\varphi \in \F(P,S)$ we obtain a morphism $r_{\varphi} : A(S) \to A(P)$.

From now on, we will use $X$ to denote a general element of $A(S)$, i.e. a linear combination of isomorphism classes of $S$-sets.
\begin{defn}
An element $X \in A(S)$ is said to be $\F$-stable if for every subgroup $P\leq S$ and every morphism $\varphi \in \F(P,S)$, we have $$ r_{\varphi} (X) = r_{i_{P,S}}(X) \quad \ \textup{in } A(P),$$ where $i_{P,S}$ is the inclusion $P \hookrightarrow S$.
\end{defn}

The following characterisation will be helpful. For the proof, we refer the reader to \cite{Reeh1}.

\begin{prop}[\cite{Reeh1}, Lemma 4.1] \label{F-stable}
Let $X$ be an element of $A(S)$. The following are equivalent:
\begin{enumerate}[$(1)$]
\item $X$ is $\F$-stable.
\item $\Phi_P(X) = \Phi_Q(X)$ whenever $P$ and $Q$ are $\F$-isomorphic subgroups of $S$.
\end{enumerate}
\end{prop}

As a corollary, we can define the Burnside ring of $\F$.

\begin{defn}
The set of $\F$-stable elements is a subring of $A(S)$, which we call the Burnside ring of $\F$ and denote by $A(\F)$.
\end{defn}

Let $\F$ be a saturated fusion system over a finite $p$-group $S$. In this case, S. P. Reeh proved that there exists a $\Z$-linear basis of $A(\F)$ formed with irreducible $\F$-stable elements of $A_+(S)$, ``irreducible'' meaning here that we cannot decompose any of these basis elements as a non-trivial sum of two $\F$-stable elements of $A_+(S)$. The following proposition will be both our main tool and the reason why we assume our fusion systems to be saturated in the remainder of the paper.

\begin{prop} \label{base alpha}
Let $\F$ be a saturated fusion system over a finite $p$-group $S$. For each $[P]_\F \in Cl(\F)$, there exists an $\F$-stable element $\alpha_{[P]_\F} \in A_+(S)$ such that the following properties are satisfied:
\begin{enumerate}[$(i)$]
\item $\alpha_{[P]_\F}$ is a linear combination of the $[S/P']$ for $P' \lesssim_\F P$.
\item There is \footnotemark a subgroup $P_0 \in [P]_\F$ such that $\alpha_{[P]_\F}$ has exactly one transitive component of the form $[S/P_0]$, while $[S/P_0]$ does not appear as a transitive component of $\alpha_{[Q]_\F}$ when $[Q]_\F \neq [P]_\F$.
\item $\Phi_{P_0} \big(\alpha_{[P]_\F} \big) = \frac{|N_S(P_0)|}{|P_0|}$.
\item $\Phi_Q \big( \alpha_{[P]_\F} \big) \neq 0$ if and only if $Q \lesssim_\F P$.
\end{enumerate}
The family $\cB := \{ \alpha_{[P]_\F} \ \vert \ [P]_\F \in Cl(\F) \}$ is a $\Z$-linear basis of $A(\F)$.
\end{prop}

\footnotetext{In fact, this property and the following are satisfied for any fully normalized subgroup in $[P]_\F$, i.e. any $P_0 \in [P]_\F$ such that $|N_S(P_0)|$ is maximal among $\F$-conjugates of $P$.}

\begin{dem}
According to Proposition 4.8 in \cite{Reeh1}, there exists, for each $[P]_\F \in Cl(\F)$, an $\F$-stable element $\alpha_{[P]_\F} \in A_+(S)$ such that 
\begin{enumerate}[$(1)$]
\item $\Phi_Q \big( \alpha_{[P]_\F} \big) = 0$ unless $Q \lesssim_\F P$.
\item For all $P'$ such that $|N_S(P')|$ is maximal in $[P]_\F$, $\alpha_{[P]_\F}$ has exactly one transitive component of the form $[S/P']$, while $\Phi_{P'} \big(\alpha_{[P]_\F} \big) = \frac{|N_S(P')|}{|P'|}$.
\item For all $P'$ such that $|N_S(P')|$ is maximal in $[P]_\F$, $[S/P']$ does not appear as a transitive component of $\alpha_{[Q]_\F}$ when $[Q]_\F \neq [P]_\F$.
\end{enumerate}
Point $(1)$ is in fact an equivalence. Indeed, fix $[P]_\F \in Cl(\F)$ and take a representative $P$ of $[P]_\F$ such that $|N_S(P)|$ is maximal. Now if $Q \lesssim_\F P$ then there exists $Q' \sim_\F Q$ such that $Q' \leq P$, and thus we get: \begin{align*}
\Phi_Q \big(\alpha_{[P]_\F} \big) &= \Phi_{Q'} \big(\alpha_{[P]_\F} \big) \qquad \qquad \, \textup{since $\alpha_{[P]_\F}$ is $\F$-stable} \\
&= \Phi_{Q'}([S/P]+X) \qquad \textup{for a certain $X \in A_+(S)$, because of point $(2)$} \\
&\geq \Phi_{Q'}([S/P]) \qquad \qquad \, \textup{because $X$ is the isomorphism class of an $S$-set} \\
&> 0 \qquad \qquad \qquad \qquad \ \ \textup{since } P' \leq Q \ .
\end{align*}

Moreover, the fact that $\cB = \{ \alpha_{[P]_\F} \ \vert \ [P]_\F \in Cl(\F) \}$ is a $\Z$-linear basis of $A(\F)$ is Corollary 4.11 in the same paper. Thus there only remains to prove Property $(i)$, which -- though it was not stated clearly in \cite{Reeh1} -- was immediate in the procedure of construction of elements $\alpha_{[P]_\F}$. It is also an easy consequence of Property $(iv)$, so we give it a proof.

Fix $[P]_\F \in Cl(\F)$, and take a transitive component $[S/Q]$ of $\alpha_{[P]_\F}$, with $Q \leq S$. Then we have $\Phi_Q \big(\alpha_{[P]_\F} \big) \geq \Phi_Q \big([S/Q] \big)$ because $\alpha_{[P]_\F}$ is the isomorphism class of an $S$-set, and thus $\Phi_Q \big(\alpha_{[P]_\F} \big) >0$ because of Lemma \ref{formule phi}. Thanks to Property $(iv)$, we deduce that $Q \lesssim_\F P$, and this concludes the proof.
\end{dem}

\subsection{The $p$-localized Burnside ring}

For various reasons, when dealing with the Burnside ring of a saturated fusion system $\F$ on a $p$-group $S$, it is often useful to work in the $p$-localisation of the $\Z$-module $A(\F)$. Thus we define the $p$-localized Burnside ring of $\F$ to be the following $\Z_{(p)}$-algebra: $$ A(\F)_{(p)} := \Z_{(p)} \otimes_\Z A(\F) $$
The Burnside ring $A(\F)$ is embedded in $A(\F)_{(p)}$ via $X \mapsto 1\otimes X$, and $A(\F)_{(p)}$ is a free $\Z_{(p)}$-module with basis $\cB = \{ \alpha_{[P]_\F} \ \vert \ [P]_\F \in Cl(\F) \}$.

In particular for $\F = \F_S(S)$, the $p$-localized Burnside ring $A(\F_S(S))_{(p)}$ is the $p$-localisation of $A(S)$, i.e. $A(S)_{(p)}$. The mark homomorphism $\Phi : A(S) \to \prod\limits_{[P]_S \in Cl(S) } \Z $ naturally extends to a homomorphism of $\Z_{(p)}$-algebras $A(S)_{(p)} \to \prod\limits_{[P]_S \in Cl(S) } \Z_{(p)}$, which we still denote by $\Phi$.

\section{Description of the prime ideals of $A(\F)$}
\label{sec descri}

In the following, $S$ is a finite $p$-group and $\F$ is a saturated fusion system over $S$. The elements $\alpha_{[P]_\F}$ for $[P]_\F \in Cl(\F)$, are the ones mentioned in Proposition \ref{base alpha}. We order them with respect to $\F$-subconjugation of the indexing subgroups, so that for each $[P]_\F,[Q]_\F \in Cl(\F)$ we have: $$ \alpha_{[P]_\F} \leq \alpha_{[Q]_\F} \iff [P]_\F \lesssim [Q]_\F \iff P \lesssim_\F Q \iff \F(P,Q) \neq \emptyset $$ 

\begin{Not}
For each $P \leq S$ and any prime $q$, we denote by $\cI_{P,q}$ and $\cI_{P,0}$ the following subsets of $A(\F)$: $$ \cI_{P,q} := \{ X \in A(\F) \ \vert \ \Phi_P(X) \equiv 0 \mod q \} $$
$$\cI_{P,0} := \{ X \in A(\F) \ \vert \ \Phi_P(X) = 0 \} $$
\end{Not}

\begin{lem}
Let $P$ be a subgroup of $S$ and $q$ be $0$ or a prime number. Then $\cI_{P,q}$ is a prime ideal of $A(\F)$. When $q \neq 0$, we actually have an isomorphism of rings: $$ A(\F)/\cI_{P,q} \simeq \Z/q\Z$$
\end{lem}

\begin{dem}
The morphism $\Phi_P$ is a ring homomorphism. If we compose it with the projection $\Z \to \Z/q\Z$ (i.e. with $\Id_\Z$ when $q=0$), we obtain a ring homomorphism $$\overline{\Phi}_P : A(\F) \longrightarrow \Z/q\Z$$ whose kernel is $\cI_{P,q}$ by definition. Hence $\cI_{P,q}$ is an ideal of $A(\F)$.

When $q$ is 0 or a prime number, $\Z/q\Z$ is an integral domain, so $\cI_{P,q}$ is a prime ideal. Moreover the image of $\overline{\Phi}_P$ is a subring of $\Z/q\Z$. When $q$ is a prime number, this image is necessarily equal to $\Z/q\Z$, and we obtain the expected isomorphism by applying the isomorphism theorem to $\overline{\Phi}_P$.
\end{dem}

The preceding lemma gives us a series of examples of prime ideals in $A(\F)$, which are the ideals $\cI_{P,q}$. In the following proposition, we show that every prime ideal of $A(\F)$ is in fact of this form. For the proof we modestly borrow our arguments from A. Dress, adapting his proof of Proposition 1. (a) in \cite{Dress}, in which he described the same result as ours, but for the Burnside ring of a finite group.

\begin{prop} \label{idéaux}
Let $\cI$ be a prime ideal of $A(\F)$. Let us denote $\cB = \{\alpha_{[P]_\F} \ \vert \ [P]_\F \in Cl(\F) \}$. Then $\cB \setminus (\cB \cap \cI)$ contains a unique minimal element $\alpha_0$, and for all $P \leq S$ such that $\alpha_{[P]_\F}=\alpha_0$, we have $$ \cI = \cI_{P,q} \ , \quad \textup{where} \quad q = \car(A(\F)/\cI) \ .$$
\end{prop}

\begin{dem}
Since $\cI$ is a prime ideal, we have $\cI \neq A(\F)$, hence $\cB \setminus (\cB \cap \cI)$ is nonempty. Assume that $\alpha_{[P]_\F}$ and $\alpha_{[Q]_\F}$ are two minimal elements of $\cB \setminus (\cB \cap \cI)$. As $\alpha_{[P]_\F}$ is a linear combination of the $[S/P']$ for $P' \lesssim_\F P$ (Property $(i)$ in Proposition \ref{base alpha}), we have: \begin{align*}
\alpha_{[P]_\F}\times \alpha_{[Q]_\F} &= \left(\sum\limits_{[P']_S \lesssim [P]_S } a_{[P']_S} \cdot [S/P'] \right) \times \left(\sum\limits_{[Q']_S \lesssim [Q]_S } b_{[Q']_S} \cdot [S/Q'] \right)\\
&= \sum\limits_{[P']_S \lesssim [P]_S } \sum\limits_{[Q']_S \lesssim [Q]_S } a_{[P']_S} b_{[Q']_S} \cdot \big([S/P'] \times [S/Q']\big) \\
& = \sum\limits_{[R]_S \lesssim [P]_S \atop [R]_S \lesssim [Q]_S } c_{[R]_S} \cdot [S/R] \qquad \qquad \textup{according to the product formula (Lemma \ref{product formula}).}
\end{align*}
So $\alpha_{[P]_\F}\times \alpha_{[Q]_\F}$ is an element of $A(\F)$ which is a linear combination of the transitive components $[S/R]$ such that $[R]_S \lesssim [P]_S$ and $[R]_S \lesssim [Q]_S$. Because of Property $(ii)$ in Proposition \ref{base alpha}, this implies that it is a linear combination of the elements $\alpha_{[R]_\F}$ with $[R]_\F \lesssim [P]_\F$ and $[R]_\F \lesssim [Q]_\F$.

Additionally we have $\alpha_{[P]_\F} \times \alpha_{[Q]_\F} \notin \cI$ because $\cI$ is a prime ideal. We deduce that there exists a class $[R]_\F \in Cl(\F)$ such that $[R]_\F \lesssim [P]_\F$, $[R]_\F \lesssim [Q]_\F$ and $\alpha_{[R]_\F} \notin \cI$. By minimality this implies $[P]_\F=[R]_\F=[Q]_\F$. This concludes the proof for unicity.

Now let us call $\alpha_0$ the minimal element of $\cB \setminus (\cB \cap \cI)$. Take $X \in A(\F)$ and fix a $P_0 \leq S$ such that $\alpha_{[P_0]_\F}=\alpha_0$ and $\alpha_0$ has exactly one transitive component equal to $[S/P_0]$ (see Property $(ii)$ in Proposition \ref{base alpha}). Since $\alpha_0$ has only transitive components of the form $[S/P']$ with $P' \lesssim_\F P_0$, we can apply a similar reasoning as above to get $$ \alpha_0 \times X = \mu \cdot \alpha_0 + \sum\limits_{[Q]_\F \lesssim [P_0]_\F \atop [Q]_\F \neq [P_0]_\F} \lambda_{[Q]_\F} \cdot \alpha_{[Q]_\F} $$ for certain $\lambda_{[Q]_\F} \in \Z$ and $\mu \in \Z$.

On the one hand $\Phi_{P_0}(\alpha_0 \times X) = \Phi_{P_0}(\alpha_0) \cdot \Phi_{P_0}(X)$ because $\Phi_{P_0}$ is a ring homomorphism, and on the other hand \begin{align*}
\Phi_{P_0}(\alpha_0 \times X) &= \Phi_{P_0}\left(\mu \cdot \alpha_0 + \sum\limits_{[Q]_\F \lesssim [P_0]_\F \atop [Q]_\F \neq [P_0]_\F} \lambda_{[Q]_\F} \cdot  \alpha_{[Q]_\F} \right) \\
&= \mu \cdot \Phi_{P_0}(\alpha_0) + \sum\limits_{[Q]_\F \lesssim [P_0]_\F \atop [Q]_\F \neq [P_0]_\F} \lambda_{[Q]_\F} \cdot  \Phi_{P_0} \big( \alpha_{[Q]_\F} \big) \qquad \textup{by linearity of }\Phi_{P_0} \\
&= \mu \cdot \Phi_{P_0}(\alpha_0)
\end{align*}
where the last equality holds because of Property $(iv)$ in Proposition \ref{base alpha}. 

Moreover, $\Phi_{P_0}(\alpha_0) \neq 0$ by Property $(iii)$ in Proposition \ref{base alpha}, so $\mu = \Phi_{P_0}(X)$. Then we obtain \begin{align*}
\alpha_0 \times X &= \Phi_{P_0}(X) \cdot \alpha_0 + \sum\limits_{[Q]_\F \lesssim [P_0]_\F \atop [Q]_\F \neq [P_0]_\F} \lambda_{[Q]_\F} \cdot \alpha_{[Q]_\F} \\
&\equiv \Phi_{P_0}(X) \cdot \alpha_0 \mod \cI  \qquad \qquad \textup{by minimality of } \alpha_0.
\end{align*}
Finally, \begin{align*}
X \in \cI &\iff \alpha_0 \times X \in \cI \qquad \qquad \textup{because } \alpha_0 \notin \cI \\
&\iff \Phi_{P_0}(X) \cdot \alpha_0 \in \cI \\
&\iff \Phi_{P_0}(X) \cdot \overline{\alpha_0} = \overline{0} \in A(\F)/\cI \\
& \iff \Phi_{P_0}(X) \equiv 0 \mod q \qquad \quad \textup{with } \ q = \car(A(\F)/\cI) \ ,
\end{align*}
following the convention that $a \equiv b \mod 0 \ $ means $a=b$ in $\Z$.
\end{dem}

Now we know that the complete collection of prime ideals of $A(\F)$ is given by the family $\cI_{P,q}$ when $q$ is prime or zero and $P \leq S$. 
In addition, thanks to Proposition \ref{F-stable} we know that $\cI_{P,q} = \cI_{Q,q}$ whenever $P \sim_\F Q$, and the previous Proposition \ref{idéaux} implies that $\cI_{P,q} \neq \cI_{P,q'}$ whenever $q \neq q'$, because the characteristic of a ring is unique. However, there could be different classes $[P]_\F$ and $[Q]_\F$ such that $\cI_{P,q} = \cI_{Q,q}$.

Before we go further, let us highlight two particular cases.

\subsection*{Particular case $q=p$}

Here we focus on prime ideals of the form $\cI_{P,p}$. For any two subgroups $P,Q \leq S$, we have $$\Phi_P([S/Q]) = \frac{|N_S(P,Q)|}{|Q|} $$ according to Lemma \ref{formule phi}. Moreover, there is a left action of $N_S(Q)$ on $N_S(P,Q)$ by multiplication, and this action is free, so $$\Phi_P([S/Q]) = \frac{|N_S(Q)|}{|Q|} \cdot k \qquad \textup{for a certain } k \in \N . $$
In particular for $Q \neq S$ we have $$\Phi_P([S/Q]) \equiv 0 \mod p \ .$$ But $\alpha_{[S]_\F} = [S/S]$ is the only isomorphism class of $S$-set in $\cB$ which contains a one-point orbit $[S/S]$ according to Property $(ii)$ in Proposition \ref{base alpha}. Hence there exists only one "type $p$" ideal, i.e. associated to the prime $p$, that is:
$$ \forall P \leq S, \qquad \cI_{P,p} = \cI_{S,p} $$ 
Since $\alpha_{[S]_\F}$ is the only element of $\cB \setminus (\cB \cap \cI_{P,p})$, the element $\alpha_{[S]_\F}$ is in particular the $\alpha_0$ in Proposition \ref{idéaux}.

Finally, we have already seen that $\{\alpha_{[P]_\F} \ \vert \ P<S \} \subseteq \cI_{S,p}$, and we clearly have $p \cdot \alpha_{[S]_\F} \in \cI_{S,p}$. The ideal generated by $\{ \alpha_{[P]_\F} \ \vert \ P<S \}\cup \{p \cdot \alpha_{[S]_\F} \}$ being clearly maximal, we deduce that $$ \cI_{S,p} = \left\langle \{\alpha_{[P]_\F} \ \vert \ P<S \}\cup \{p \cdot \alpha_{[S]_\F} \} \right\rangle \ .$$
The ideal $\cI_{S,p}$ is in particular maximal.

\subsection*{Particular case $q=0$}

The other case which stands out is $\cI_{P,0}$. Let $P$ be a subgroup of $S$. According to Property $(iv)$ in Proposition \ref{base alpha}, we know that for every $Q \leq S$, $$ \alpha_{[Q]_\F} \notin \cI_{P,0} \iff \Phi_P\big(\alpha_{[Q]_\F} \big) \neq 0 \iff P \lesssim_\F Q \ .$$

From there, we deduce that $$\cB \setminus (\cB \cap \cI_{P,0}) = \{ \alpha_{[Q]_\F} \ \vert \ P \lesssim_\F Q \} \ .$$
We also see that each class $[P]_\F \in Cl(\F)$ gives rise to a distinct $\cI_{P,0}$.

As a consequence of the above, we have the following inclusion of ideals $$ \left\langle \alpha_{[Q]_\F} \ \vert \ P \text{ is not $\F$-subconjugate to } Q \right\rangle \subseteq \cI_{P,0} \ .$$ This inclusion is in fact strict for each $P<S$, because in such a case $\cI_{P,0}$ also contains elements that do not belong to the left-hand ideal, such as $ \alpha_{[P]_\F} - \frac{|N_S(P_0)|}{|P_0|} \cdot \alpha_{[S]_\F} $ with $P_0$ fully normalized in $[P]_\F$ (see footnote in Proposition \ref{base alpha}). Nonetheless, we can exhibit a basis for $\cI_{P,0}$. Indeed, for any $X \in A(\F)$, one can write $$ X \ = \ \sum_{[Q]_{\F} \in Cl(\F)} \lambda_{[Q]_{\F}} \cdot \alpha_{[Q]_\F} \qquad \textup{with} \quad \lambda_{[Q]_{\F}} \in \Z \ .$$
By linearity of $\Phi_P$, and because $\Phi_P\big( \alpha_{[Q]_\F} \big) = 0$ for any $Q$ such that $P \not\lesssim_\F Q$, the property of belonging to $\cI_{P,0}$ is equivalent to the equation: $$ \sum_{[P]_\F \lesssim [Q]_{\F}} \Phi_P\big( \alpha_{[Q]_\F} \big) \cdot \lambda_{[Q]_{\F}} \ = \ 0 $$
Now, since $\Phi_P \big(\alpha_{[S]_\F} \big) = 1$, the equation rephrases as $$ \lambda_{[S]_\F} \ = \ - \sum_{\substack{[P]_\F \lesssim [Q]_{\F} \\ [Q]_\F \neq [S]_\F}} \Phi_P\big( \alpha_{[Q]_\F} \big) \cdot \lambda_{[Q]_{\F}} \ , $$ so $X \in \cI_{P,0}$ can be expressed as follows:
$$ X \ = \ \sum_{[P]_\F \not\lesssim [Q]_{\F}} \lambda_{[Q]_{\F}} \cdot \alpha_{[Q]_\F} \ + \ \sum_{\substack{[P]_\F \lesssim [Q]_{\F} \\ [Q]_\F \neq [S]_\F}} \lambda_{[Q]_\F} \cdot \Big( \alpha_{[Q]_\F} - \Phi_P\big( \alpha_{[Q]_\F} \big) \cdot \alpha_{[S]_\F} \Big) $$

Moreover, as $\{\alpha_{[Q]_\F} \ \vert \ [Q]_\F \in Cl(\F) \}$ is a basis of $A(\F)$, we know that $X$ equals $0$ only if $\lambda_{[Q]_{\F}} =0$ for any $[Q]_\F$. Hence the family $$ \left\{ \alpha_{[Q]_\F} \ \Big| \ [P]_\F \not\lesssim [Q]_\F \right\} \cup \left\{ \alpha_{[Q]_\F} - \Phi_P\big( \alpha_{[Q]_\F} \big) \cdot \alpha_{[S]_\F} \ \Big| \ [P]_\F \lesssim [Q]_\F \ \textup{and} \ [Q]_\F \neq [S]_\F \right\} $$ is a basis of $\cI_{P,0}$. Keeping in mind that $\Phi_P\big( \alpha_{[Q]_\F} \big) = 0$ for any $Q$ such that $P \not\lesssim_\F Q$, we eventually have: \begin{align*}
\cI_{P,0} \ &= \ \Big\langle \left\{ \alpha_{[Q]_\F} \ \Big| \ [P]_\F \not\lesssim [Q]_\F \right\} \cup \left\{ \alpha_{[Q]_\F} - \Phi_P\big( \alpha_{[Q]_\F} \big) \cdot \alpha_{[S]_\F} \ \Big| \ [P]_\F \lesssim [Q]_\F \ \textup{and} \ [Q]_\F \neq [S]_\F \right\} \Big\rangle \\[0.5em]
&= \ \Big\langle \left\{ \alpha_{[Q]_\F} - \Phi_P\big( \alpha_{[Q]_\F} \big) \cdot \alpha_{[S]_\F} \ \Big| \ [Q]_\F \neq [S]_\F \right\} \Big\rangle \\
&= \ \Big\langle \left\{ \alpha_{[Q]_\F} - \Phi_P\big( \alpha_{[Q]_\F} \big) \cdot \alpha_{[S]_\F} \ \Big| \ Q<S \right\} \Big\rangle
\end{align*}

Notice that, for $P=S$, we get: $$ \cI_{S,0} = \left\langle \alpha_{[Q]_\F} \ \vert \ Q < S \right\rangle $$

\subsection*{Exhaustive description}

Now we are ready to refine our description of the prime ideals of $A(\F)$.

\begin{thm} \label{typology}
Let $n$ denote the cardinality of $Cl(\F)$. The prime ideals of $A(\F)$ are the following:
\begin{enumerate}[$\bullet$]
\item \textup{\textbf{"Type $p$" ideal.}} There is only one "type $p$" prime ideal, which is $$\cI_{S,p} \ = \ \{ X \in A(\F) \ \vert \ \Phi_S(X) \equiv 0 \mod p \} \ = \ \left\langle \{\alpha_{[P]_\F} \ \vert \ P<S \}\cup \{p \cdot \alpha_{[S]_\F} \} \right\rangle \ .$$
\item \textup{\textbf{"Type $0$" ideals.}} There are exactly $n$ "type $0$" prime ideals, one per class $[P]_\F \in Cl(\F)$: $$ \cI_{P,0} \ = \ \{ X \in A(\F) \ \vert \ \Phi_P(X) =0 \} \ = \ \left\langle \left\{ \alpha_{[Q]_\F} - \Phi_P\big( \alpha_{[Q]_\F} \big) \cdot \alpha_{[S]_\F} \ \big| \ Q<S \right\} \right\rangle \ .$$
\item \textup{\textbf{"Type $q$" ideals, for $q \neq p$ prime.}} There are exactly $n$ "type $q$" prime ideals, one per class $[P]_\F$: $$ \cI_{P,q} \ = \ \{ X \in A(\F) \ \vert \ \Phi_P(X) \equiv 0 \mod q \} \ .$$
\end{enumerate}
\end{thm}

\begin{dem}
According to Proposition \ref{idéaux}, we know that every prime ideal is of the form $\cI_{P,q}$ for some $P \leq S$ and some $q$ prime or zero. We call $q$ the type of the ideal. We also know that different types give different prime ideals, because the type is the characteristic of $A(\F)/\cI_{P,q}$, which is unique.

Moreover, we already proved the classification of the types $q=p$ and $q=0$ in our study of these particular cases. Thus what remains to prove is the general case of a prime $q \neq p$. Let $P$ and $Q$ be subgroups of $S$. We want to prove that $$ P \sim_\F Q \iff \cI_{P,q} = \cI_{Q,q} \ .$$
The left-to-right implication is a consequence of Proposition \ref{F-stable}. For the other direction, consider the elements $\alpha_{[P]_\F}$ and $\alpha_{[Q]_\F}$. Up to replacing $P$ and $Q$ by $\F$-isomorphic subgroups, we may assume that $$ \Phi_P \big(\alpha_{[P]_\F} \big) = \frac{|N_S(P)|}{|P|} \qquad \textup{and} \qquad \Phi_Q \big(\alpha_{[Q]_\F} \big) = \frac{|N_S(Q)|}{|Q|} $$ according to Properties $(ii)$ and $(iii)$ in Proposition \ref{base alpha}.

As a consequence, since $S$ is a $p$-group and $q \neq p$, both $\Phi_P \big(\alpha_{[P]_\F} \big)$ and $\Phi_Q \big(\alpha_{[Q]_\F} \big)$ are non-zero in $\Z/q\Z$. Hence $$ \alpha_{[P]_\F} \notin \cI_{P,q} \qquad \textup{and} \qquad \alpha_{[Q]_\F} \notin \cI_{Q,q} \ , $$ and because of the equality $\cI_{P,q} = \cI_{Q,q}$ we deduce that $$ \alpha_{[P]_\F} \notin \cI_{Q,q} \qquad \textup{and} \qquad \alpha_{[Q]_\F} \notin \cI_{P,q} \ ,$$ or in other words $$ \Phi_Q \big(\alpha_{[P]_\F} \big) \not\equiv 0 \mod q \qquad \textup{and} \qquad \Phi_P \big(\alpha_{[Q]_\F} \big) \not\equiv 0 \mod q \ . $$
Then Property $(iv)$ in Proposition \ref{base alpha} implies that $$ Q \lesssim_\F P \qquad \textup{and} \qquad P \lesssim_\F Q \ , \qquad \textup{i.e.} \qquad P \sim_\F Q \ .$$
\end{dem}

\begin{Not}
As a consequence of Theorem \ref{typology}, and in order to have a unique notation for each ideal, we will exclusively use the notation $\cI_{S,p}$ to nominate the unique "type $p$" ideal, and no more $\cI_{P,p}$ when $P < S$. Thus, if we denote by $\cP_{+0}^{-p}$ the set of prime numbers from which we have removed $p$ and added 0, and if we pick a set $\cX_\F$ of representatives of classes in $Cl(\F)$, we get a bijection between the set of indices $$ \big( \cX_\F \times \cP_{+0}^{-p} \big) \cup \big\{(S,p)\big\}$$
and the set of prime ideals of $A(\F)$, given by $(P,q) \mapsto \cI_{P,q}$.
\end{Not}

\section{Inclusion between prime ideals of $A(\F)$}
\label{sec inclusions}

In this section we look at the inclusion relations between prime ideals of $A(\F)$. It is obvious that we always have the inclusion $\cI_{P,0} \subseteq \cI_{P,q}$ (or $\cI_{P,0} \subseteq \cI_{S,p}$ for the case $q=p$). Whenever $q \neq 0$ this inclusion is strict, because $\cI_{P,q}$ also contains (among others) all elements of the form $q \cdot \alpha_{[Q]_\F}$, especially when $P \lesssim_\F Q$.

\begin{prop} \label{inclusions}
The following equivalence describes all possible inclusions between the ideals $\cI_{P,q}$, recalling that whenever $q=p$, with this notation, we have $P=S$.
$$ \cI_{P,q} \subseteq \cI_{P',q'} {\bm \iff} (q=q' \textup{ and } P \sim_\F P') \textup{ \textbf{or} } (q=0, \  q'\notin \{0,p\} \textup{ and } P \sim_\F P') \textup{ \textbf{or} } (q=0 \textup{ and } q'=p) $$
In particular, $\cI_{P,q}$ is a minimal (resp. maximal) prime ideal if and only if $q=0$ (resp. $q \neq 0$).
\end{prop}

\begin{dem}
The right-to-left direction of the equivalence has already been discussed. For the other direction, take two prime ideals $\cI_{P,q} \subseteq \cI_{P',q'}$. This inclusion implies that we have a surjective ring homomorphism $$ A(\F)/\cI_{P,q} \longrightarrow A(\F)/\cI_{P',q'} \ .$$ Because of the characteristic of those rings, this is possible only if $q'$ divides $q$, i.e. $q=0$ or $q'=q$. The case $q'=q=p$ is trivial, because then $P'=P=S$ so in particular $P \sim_\F P'$. Hence there only remains to show that whenever $q' \neq p$ we have $P \sim_\F P'$.

Assume that $q' \neq p$. We introduce the following element: $$ X=\Phi_P \big(\alpha_{[P]_\F} \big) \cdot \alpha_{[P']_\F} - \Phi_P \big(\alpha_{[P']_\F} \big) \cdot \alpha_{[P]_\F} \in A(\F) \ . $$ We have $$ \Phi_P(X) = \Phi_P \big(\alpha_{[P]_\F}\big) \cdot \Phi_P \big(\alpha_{[P']_\F}\big) - \Phi_P \big(\alpha_{[P']_\F}\big) \cdot \Phi_P \big(\alpha_{[P]_\F} \big) = 0 \ ,$$ 
so $X$ is an element of $\cI_{P,0}$, which is included in $\cI_{P,q}$. Since we also have $\cI_{P,q} \subseteq \cI_{P',q'}$, it means that $\Phi_{P'}(X) \equiv 0 \mod q'$. On the other hand, $$ \Phi_{P'}(X) =  \Phi_P \big(\alpha_{[P]_\F} \big) \cdot \Phi_{P'} \big(\alpha_{[P']_\F}\big) - \Phi_P \big(\alpha_{[P']_\F}\big) \cdot \Phi_{P'} \big(\alpha_{[P]_\F} \big) \ .$$ Up to replacing $P$ and $P'$ by $\F$-isomorphic subgroups, thanks to Proposition \ref{base alpha} we can assume that $$ \Phi_P \big(\alpha_{[P]_\F} \big) = \frac{|N_S(P)|}{|P|} \qquad \textup{and} \qquad \Phi_{P'} \big(\alpha_{[P']_\F} \big) = \frac{|N_S(P')|}{|P'|} \ .$$ As $q'\neq p$, this gives $$ \Phi_P \big(\alpha_{[P]_\F} \big) \cdot \Phi_{P'} \big(\alpha_{[P']_\F}\big) \ \not\equiv \ 0 \mod q' $$ and so $$ \Phi_P \big(\alpha_{[P']_\F}\big) \cdot \Phi_{P'} \big(\alpha_{[P]_\F} \big) = \Phi_P \big(\alpha_{[P]_\F} \big) \cdot \Phi_{P'} \big(\alpha_{[P']_\F}\big) - \Phi_{P'}(X) \ \not\equiv \ 0 \mod q' \ .$$
In particular, both $\Phi_P \big(\alpha_{[P']_\F}\big)$ and $\Phi_{P'} \big(\alpha_{[P]_\F} \big)$ have to be non-zero, which implies $P \sim_\F P'$ according to Property $(iv)$ in Proposition \ref{base alpha}.
\end{dem}

Notice that, in the above proposition, all the inclusions $\cI_{P,q} \subseteq \cI_{P',q'}$ are strict except for the case when $q=q'$ and $P \sim_\F P'$.

\section{Prime ideals of $A(\F)_{(p)}$}
\label{sec local}

\subsection{Classification}

In this section we are interested in the prime ideals of the $p$-localized Burnside ring $A(\F)_{(p)}$. It appears that the proof of our Proposition \ref{idéaux} is still applicable in the $p$-localized case because it essentially relies on the properties of the basis $\cB$, which is a basis for both $A(\F)$ and $A(\F)_{(p)}$ (as a $\Z$-module and as a $\Z_{(p)}$-module respectively). 

Yet there is a noticeable difference in the $p$-localized case: for each prime $q$ different from $p$, $q$ becomes invertible in $A(\F)_{(p)}$, so that the only possible characteristics for a quotient of $A(\F)_{(p)}$ are $0$ and $p$. Consequently, the prime ideals of $A(\F)_{(p)}$ are necessarily "type $0$" ideals or "type $p$" ideals.

\begin{prop} \label{idéaux localisés}
Let $\cJ$ be a prime ideal of $A(\F)_{(p)}$. Then $\cB \setminus (\cB \cap \cJ)$ contains a unique minimal element $\alpha_0$, and there are two possible cases, depending on the type $q := \car(A(\F)_{(p)}/\cJ)$:
\begin{enumerate}[\ding{250}]
\item Either $q=p$, then $\alpha_0 = \alpha_{[S]_\F}$ and we have: $$ \cJ = \{ X \in A(\F)_{(p)} \ \vert \ \Phi_S(X) \equiv 0 \mod p \} \ =: \ \cJ_{S,p} \ . $$
\item Or $q=0$, and for each $P \leq S$ such that $\alpha_{[P]_\F} = \alpha_0$ we have: $$ \cJ = \{ X \in A(\F)_{(p)} \ \vert \ \Phi_P(X) =0 \} \ =: \ \cJ_{P,0} \ . $$
\end{enumerate}
\end{prop}

\begin{dem}
As we said above, the proof of Proposition \ref{idéaux} still applies. However we can give another proof, using the properties of $p$-localisation. In fact, as $A(\F)_{(p)}$ is the localisation of $A(\F)$ by the multiplicative set $$ \cT := \{ n \cdot [S/S] \ \vert \ n \in \Z \setminus p\Z \} \ ,$$ there is a bijection between prime ideals of $A(\F)_{(p)}$ and prime ideals of $A(\F)$ that do not intersect $\cT$. This bijection is given by $\cJ \mapsto \cJ \cap A(\F)$ for the direct sense and $\cI \mapsto \cI_{(p)}$ for the converse.

Now if we look at our classification of prime ideals of $A(\F)$ in Theorem \ref{typology}, we remark that the ideals which do not intersect $\cT$ are precisely those with type $0$ or $p$. Because the mark homomorphism commutes with $p$-localisation, those ideals are: $$ \big(\cI_{S,p}\big)_{(p)} = \{ X \in A(\F)_{(p)} \ \vert \ \Phi_S(X) \equiv 0 \mod p \} = \cJ_{S,p}  \,$$ and $$ \big(\cI_{P,0}\big)_{(p)} = \{ X \in A(\F)_{(p)} \ \vert \ \Phi_P(X) =0 \} = \cJ_{P,0} \ .$$
\end{dem}

The proof above also tells us that the ideals $\cJ_{P,0}$ and $\cJ_{Q,0}$ are equal if and only if $\cI_{P,0}$ and $\cI_{Q,0}$ are equal in $A(\F)$, which happens if and only if $P \sim_\F Q$. Moreover the inclusion relations between those ideals are quite obvious, so we can complete our classification of the prime ideals in $A(\F)_{(p)}$.

\begin{thm} \label{thm prime ideals p-localized}
Let $n$ denote the cardinality of $Cl(\F)$. The prime ideals of $A(\F)_{(p)}$ are the following:
\begin{enumerate}[$\bullet$]
\item \textup{\textbf{"Type $p$" ideal.}} There is only one "type $p$" prime ideal, which is maximal: $$\cJ_{S,p} \ = \ \{ X \in A(\F)_{(p)} \ \vert \ \Phi_S(X) \equiv 0 \mod p \} \ = \ \left\langle \{\alpha_{[P]_\F} \ \vert \ P<S \}\cup \{p \cdot \alpha_{[S]_\F}\} \right\rangle \ .$$
\item \textup{\textbf{"Type $0$" ideals.}} There are exactly $n$ "type $0$" prime ideals, one per class $[P]_\F \in Cl(\F)$: $$ \cJ_{P,0} \ = \ \{ X \in A(\F)_{(p)} \ \vert \ \Phi_P(X) =0 \} \ = \ \left\langle \left\{ \alpha_{[Q]_\F} - \Phi_P\big( \alpha_{[Q]_\F} \big) \cdot \alpha_{[S]_\F} \ \big| \ Q<S \right\} \right\rangle \ .$$
All "type 0" prime ideals are minimal and contained in $\cJ_{S,p}$.
\end{enumerate}
\end{thm}

\begin{dem}
The bijection, invoked in the proof of Proposition \ref{idéaux localisés}, between prime ideals of $A(\F)_{(p)}$ and prime ideals of $A(\F)$ that do not intersect $\cT$, is order-preserving for the partial order given by inclusion. Since all "type $0$" (resp. "type $p$") prime ideals in $A(\F)$ are minimal (resp. maximal) and all are distinct, the result follows for prime ideals in $A(\F)_{(p)}$.
\end{dem}

\begin{Rq}
Notice that we always have $$ \{ \alpha_{[Q]_\F} \ \vert P \lesssim_{\F} Q \} = \cB \setminus (\cB \cap \cJ_{P,0}) \qquad \textup{and} \qquad  \left\langle \alpha_{[Q]_\F} \ \vert \ P \not\lesssim_\F Q \right\rangle \subseteq \cJ_{P,0} \ . $$
\end{Rq}

\begin{cor}
The $p$-localized Burnside ring $A(\F)_{(p)}$ of a saturated fusion system is a local ring, and its residue field is isomorphic to $\bF_p$.
\end{cor}

\subsection{Examples with $S=D_8$}

As an illustration, we will describe all prime ideals of $A(\F)_{(2)}$, using generators, when $\F$ is one of the two non-inner saturated fusion systems over $D_8$.

\begin{figure}[!h]
\begin{center}
\scalebox{0.8}{
\begin{tikzpicture}[node distance=2.5cm]
\tikzset{sg/.style={draw,rectangle,align=center,fill=white}}
\tikzset{sgen/.style={draw,rectangle,drop shadow={opacity=0.8},align=center,fill=white}}
 \node[sg] (a0)                  {$\{1\}$ \\ $\langle 1 \rangle$};
 \node[sg] (a3)  [above of=a0]   {$Z \simeq \Z/2\Z$ \\ $\langle r^2 \rangle$};
 \node[sg] (a2)  [left of=a3]   {$\Z/2\Z$ \\ $\langle rs \rangle$};
 \node[sg] (a1)  [left of=a2]   {$\Z/2\Z$ \\ $\langle r^{-1}s \rangle$};
 \node[sg] (a4)  [right of=a3]   {$\Z/2\Z$ \\ $\langle s \rangle$};
 \node[sg] (a5)  [right of=a4]   {$\Z/2\Z$ \\ $\langle r^2s \rangle$};
 \node[sg] (a7)  [above of=a3]   {$C \simeq \Z/4\Z$ \\ $\langle r \rangle$};
 \node[sg] (a6)  [left of=a7]  {$(\Z/2\Z)^2$ \\ $\langle r^2,rs \rangle$};
 \node[sg] (a8)  [right of=a7]   {$(\Z/2\Z)^2$ \\ $\langle r^2,s \rangle$};
 \node[sg] (a9)  [above of=a7]  {$D_8$ \\ $\langle r,s \rangle$};
 \draw (a0.north)   -- (a1.south);
 \draw (a0.north)   -- (a2.south);
 \draw (a0.north)   -- (a3.south);
 \draw (a0.north)   -- (a4.south);
 \draw (a0.north)   -- (a5.south);
 \draw (a1.north)   -- (a6.south);
 \draw (a2.north)   -- (a6.south);
 \draw (a3.north)   -- (a6.south);
 \draw (a3.north)   -- (a7.south);
 \draw (a3.north)   -- (a8.south);
 \draw (a4.north)   -- (a8.south);
 \draw (a5.north)   -- (a8.south);
 \draw (a6.north)   -- (a9.south);
 \draw (a7.north)   -- (a9.south);
 \draw (a8.north)   -- (a9.south);
 \draw[dashed,color=black!50,thick] (a1)   -- (a2);
 \draw[dashed,color=black!50,thick] (a4)   -- (a5);
\end{tikzpicture}
}
\caption{Subgroup lattice of $D_8$}
\label{treillis sg D_8}
\end{center}
\end{figure}
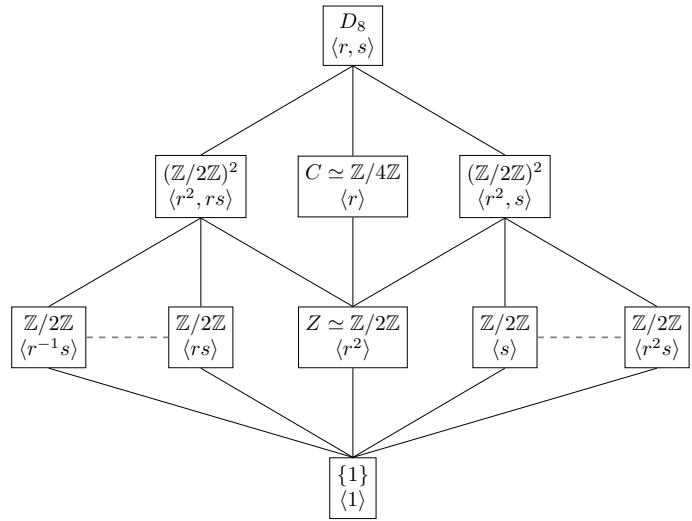

In this section, $S$ denotes the dihedral group $D_8$ of order $8$, generated by the symmetry $s$ of order $2$ and the rotation $r$ of order $4$. We also denote by $Z := \langle r^2 \rangle$ the center of $D_8$, and by $C:= \langle r \rangle$ its cyclical subgroup of order $4$. Recall that the lattice of subgroups of $D_8$ is as represented in Figure \ref{treillis sg D_8}, where the horizontal dashed lines stand for conjugation in $S$.

\noindent We choose the following order on $Cl(S)$: $$ [\triv]_S < [Z]_S < [\langle rs \rangle ]_S < [\langle s \rangle ]_S < [C]_S < [\langle r^2,rs \rangle ]_S < [\langle r^2,s \rangle ]_S < [D_8]_S $$

For each one of our two examples, we first give the expressions of the elements $\alpha_P$ using the mark homomorphism, i.e. as vectors in $\prod\limits_{[P]_S \in Cl(S)} \Z$ (with the order induced by the above order on $Cl(S)$). Those expressions are calculated using the properties in Proposition \ref{base alpha}. Then we list all prime ideals in $A(\F)_{(2)}$ and describe them with a (minimal) set of generators, expressed in terms of the $\alpha_{[P]_\F}$. The calculations lie mainly on the table of marks of the group $D_8$, i.e. the matrix $\big(\Phi_P([S/Q])\big)_{[P]_S,[Q]_S \in Cl(S)}$, so we give it below.

$$ \begin{pmatrix}
8 & 0 & 0 & 0 & 0 & 0 & 0 & 0 \\
4 & 4 & 0 & 0 & 0 & 0 & 0 & 0 \\
4 & 0 & 2 & 0 & 0 & 0 & 0 & 0 \\
4 & 0 & 0 & 2 & 0 & 0 & 0 & 0 \\
2 & 2 & 0 & 0 & 2 & 0 & 0 & 0 \\
2 & 2 & 2 & 0 & 0 & 2 & 0 & 0 \\
2 & 2 & 0 & 2 & 0 & 0 & 2 & 0 \\
1 & 1 & 1 & 1 & 1 & 1 & 1 & 1 \\
\end{pmatrix} $$

In order to enhance the readability of the description of prime ideals, we will denote $\alpha_P$ instead of $\alpha_{[P]_\F}$, which is a slight abuse of notation.

\subsubsection{Prime ideals of $A(\F)_{(2)}$ when $\F=\F_{D_8}(S_4)$}

There are two different ways to embed our $D_8$ as a $2$-Sylow subgroup of $S_4$. Here we consider the situation where the center $Z$ is $\F$-isomorphic to the subgroup $\langle rs \rangle$, i.e. the situation where $\langle r^2,rs \rangle$ is the Klein subgroup of double-transpositions in $S_4$. The values of basis elements through the mark homomorphism are given below.
\begin{align*}
\alpha_{ \triv } \quad &: \quad \begin{pmatrix} 8 & 0 & 0 & 0 & 0 & 0 & 0 & 0 \end{pmatrix} \\
\alpha_{ Z } \quad &: \quad \begin{pmatrix} 12 & 4 & 4 & 0 & 0 & 0 & 0 & 0 \end{pmatrix} \\
\alpha_{ \langle s \rangle } \quad &: \quad \begin{pmatrix} 4 & 0 & 0 & 2 & 0 & 0 & 0 & 0 \end{pmatrix} \\
\alpha_{ C } \quad &: \quad \begin{pmatrix} 6 & 2 & 2 & 0 & 2 & 0 & 0 & 0 \end{pmatrix} \\
\alpha_{ \langle r^2,rs \rangle } \quad &: \quad \begin{pmatrix} 2 & 2 & 2 & 0 & 0 & 2 & 0 & 0 \end{pmatrix} \\
\alpha_{ \langle r^2,s \rangle } \quad &: \quad \begin{pmatrix} 6 & 2 & 2 & 2 & 0 & 0 & 2 & 0 \end{pmatrix} \\
\alpha_{ D_8 } \quad &: \quad \begin{pmatrix} 1 & 1 & 1 & 1 & 1 & 1 & 1 & 1 \end{pmatrix}
\end{align*}

\noindent \emph{Maximal "type 2" ideal:}
$$ \cJ_{S,p} = \left\langle \alpha_{ \triv } \ , \ \alpha_{ Z } \ , \ \alpha_{ \langle s \rangle } \ , \ \alpha_{ C } \ , \ \alpha_{ \langle r^2,rs \rangle } \ , \ \alpha_{ \langle r^2,s \rangle } \ , \ 2\cdot \alpha_{ D_8 } \right\rangle $$
\noindent \emph{Minimal "type 0" ideals:}
\begin{align*}
\cJ_{S,0} &= \left\langle\alpha_{ \triv } \ , \ \alpha_{ Z } \ , \ \alpha_{ \langle s \rangle } \ , \ \alpha_{ C } \ , \ \alpha_{ \langle r^2,rs \rangle } \ , \ \alpha_{ \langle r^2,s \rangle } \right\rangle \\[1em]
\cJ_{\langle r^2,s \rangle,0} &= \left\langle \alpha_{ \triv } \ , \ \alpha_{ Z } \ , \ \alpha_{ \langle s \rangle } \ , \ \alpha_{ C } \ , \ \alpha_{ \langle r^2,rs \rangle } \ , \ \alpha_{ \langle r^2,s \rangle }-2 \cdot \alpha_{ D_8 } \right\rangle \\[1em]
\cJ_{\langle r^2,rs \rangle,0} &= \left\langle \alpha_{ \triv } \ , \ \alpha_{ Z } \ , \ \alpha_{ \langle s \rangle } \ , \ \alpha_{ C } \ , \ \alpha_{ \langle r^2,s \rangle } \ , \ \alpha_{ \langle r^2,rs \rangle }-2 \cdot \alpha_{ D_8 } \right\rangle \\[1em]
\cJ_{C,0} &= \left\langle \alpha_{ \triv } \ , \ \alpha_{ Z } \ , \ \alpha_{ \langle s \rangle } \ , \ \alpha_{ \langle r^2,rs \rangle } \ , \ \alpha_{  \langle r^2,s \rangle } \ , \ \alpha_{ C }-2 \cdot \alpha_{ D_8 } \right\rangle \\[1em]
\cJ_{\langle s \rangle,0} &= \left\langle \alpha_{ \triv } \ , \ \alpha_{ Z } \ , \ \alpha_{ C } \ , \ \alpha_{ \langle r^2,rs \rangle } \ , \ \alpha_{ \langle s \rangle }-2 \cdot \alpha_{ D_8 } \ , \ \alpha_{ \langle r^2,s \rangle }-2 \cdot \alpha_{ D_8 } \right\rangle \\[1em]
\cJ_{Z,0} &= \left\langle \alpha_{  \triv } \ , \ \alpha_{ \langle s \rangle } \ , \ \alpha_{ Z }-4 \cdot \alpha_{ D_8 } \ , \ \alpha_{ C }-2 \cdot \alpha_{ D_8 } \ , \ \alpha_{ \langle r^2,rs \rangle }-2 \cdot \alpha_{ D_8 } \ , \ \alpha_{ \langle r^2,s \rangle }-2 \cdot \alpha_{ D_8 } \right\rangle \\[1em]
\cJ_{\triv,0} &= \left\langle \alpha_{ \triv }-8 \cdot \alpha_{ D_8 } \ , \ \alpha_{ Z }-12 \cdot \alpha_{ D_8 } \ , \ \alpha_{ \langle s \rangle }-4 \cdot \alpha_{ D_8 } \ , \ \alpha_{ C }-6 \cdot \alpha_{ D_8 } \ , \ \alpha_{ \langle r^2,rs \rangle }-2 \cdot \alpha_{ D_8 } \ , \ \alpha_{ \langle r^2,s \rangle }-6 \cdot \alpha_{ D_8 } \right\rangle 
\end{align*}

\subsubsection{Prime ideals of $A(\F)_{(2)}$ when $\F=\F_{D_8}(A_6)$}

Here again we begin with the expression of elements $\alpha_{ P }$ through the mark homomorphism. Recall that when $D_8$ is embedded in $A_6$, all of its order $2$ subgroups are conjugate in $A_6$.
\begin{align*}
\alpha_{ \triv } \quad &: \quad \begin{pmatrix} 8 & 0 & 0 & 0 & 0 & 0 & 0 & 0 \end{pmatrix} \\
\alpha_{ Z } \quad &: \quad \begin{pmatrix} 20 & 4 & 4 & 4 & 0 & 0 & 0 & 0 \end{pmatrix} \\
\alpha_{ C } \quad &: \quad \begin{pmatrix} 10 & 2 & 2 & 2 & 2 & 0 & 0 & 0 \end{pmatrix} \\
\alpha_{ \langle r^2,rs \rangle } \quad &: \quad \begin{pmatrix} 6 & 2 & 2 & 2 & 0 & 2 & 0 & 0 \end{pmatrix} \\
\alpha_{ \langle r^2,s \rangle } \quad &: \quad \begin{pmatrix} 6 & 2 & 2 & 2 & 0 & 0 & 2 & 0 \end{pmatrix} \\
\alpha_{ D_8 } \quad &: \quad \begin{pmatrix} 1 & 1 & 1 & 1 & 1 & 1 & 1 & 1 \end{pmatrix}
\end{align*}

\noindent \emph{Maximal "type 2" ideal:}
$$ \cJ_{S,p} = \left\langle \alpha_{ \triv } \ , \ \alpha_{ Z } \ , \ \alpha_{ C } \ , \ \alpha_{ \langle r^2,rs \rangle } \ , \ \alpha_{  \langle r^2,s \rangle } \ , \ 2\cdot \alpha_{ D_8 } \right\rangle $$
\noindent \emph{Minimal "type 0" ideals:}
\begin{align*}
\cJ_{S,0} &= \left\langle \alpha_{ \triv } \ , \ \alpha_{ Z } \ , \ \alpha_{ C } \ , \ \alpha_{ \langle r^2,rs \rangle } \ , \ \alpha_{ \langle r^2,s \rangle } \right\rangle \\[1em]
\cJ_{\langle r^2,s \rangle,0} &= \left\langle \alpha_{ \triv } \ , \ \alpha_{ Z } \ , \ \alpha_{ C } \ , \ \alpha_{ \langle r^2,rs \rangle } \ , \ \alpha_{ \langle r^2,s \rangle }-2 \cdot \alpha_{ D_8 } \right\rangle \\[1em]
\cJ_{\langle r^2,rs \rangle,0} &= \left\langle \alpha_{ \triv } \ , \ \alpha_{ Z } \ , \ \alpha_{ C } \ , \ \alpha_{ \langle r^2,s \rangle } \ , \ \alpha_{ \langle r^2,rs \rangle }-2 \cdot \alpha_{ D_8 } \right\rangle \\[1em]
\cJ_{C,0} &= \left\langle \alpha_{ \triv } \ , \ \alpha_{ Z } \ , \ \alpha_{ \langle r^2,rs \rangle } \ , \ \alpha_{ \langle r^2,s \rangle } \ , \ \alpha_{ C }-2 \cdot \alpha_{ D_8 } \right\rangle \\[1em]
\cJ_{Z,0} &= \left\langle \alpha_{ \triv } \ , \ \alpha_{ Z }-4 \cdot \alpha_{ D_8 } \ , \ \alpha_{ C }-2 \cdot \alpha_{ D_8 } \ , \ \alpha_{ \langle r^2,rs \rangle }-2 \cdot \alpha_{ D_8 } \ , \ \alpha_{ \langle r^2,s \rangle }-2 \cdot \alpha_{ D_8 } \right\rangle \\[1em]
\cJ_{\triv,0} &= \left\langle \alpha_{ \triv }-8 \cdot \alpha_{ D_8 } \ , \ \alpha_{ Z }-20 \cdot \alpha_{ D_8 } \ , \ \alpha_{ C }-10 \cdot \alpha_{ D_8 } \ , \ \alpha_{ \langle r^2,rs \rangle }-6 \cdot \alpha_{ D_8 } \ , \ \alpha_{ \langle r^2,s \rangle }-6 \cdot \alpha_{ D_8 } \right\rangle 
\end{align*}

\subsubsection{About the computations}

In the above examples, the expression of the basis elements $\alpha_{[Q]_\F}$ through the mark homomorphism have been computed using the following remarks:
\begin{enumerate}[\ding{250}]
\item Each one of these vectors is a linear combination of the lines of the table of marks of $D_8$, with nonnegative integer coefficients.
\item For $\alpha_{[Q]_\F}$, these coefficients have to be zero unless they are indexed by $[P]_{S}$ with $P \lesssim_\F Q$.
\item The combination has to be minimal, so that $\alpha_{[Q]}$ is irreducible.
\item Each vector has to represent an $\F$-stable element, which means, according to Proposition \ref{F-stable}, that if two coordinates have $\F$-isomorphic indices, they have to be equal.
\end{enumerate}
Thus, to get the vector associated with $\alpha_{[Q]_\F}$, start with the line indexed by $[Q_0]_S$ in the table of marks, where $Q_0$ is a fully normalized subgroup in $[Q]_\F$ (see footnote in Proposition \ref{base alpha}), then add other lines in view of the above remarks. For the examples on $D_8$, this process is quite straightforward. For instance, if we want to compute the vector associated with $\alpha_Z$ in the case where $\F = \F_{D_8}(A_6)$, we start with the vector $$\begin{pmatrix}
4 & 4 & 0 & 0 & 0 & 0 & 0 & 0 
\end{pmatrix} $$ and we need to add other lines of the table of marks in order to get a vector with equal values on the second, the third and the fourth coordinates. The only way to do this is to add copies of the third and the fourth lines, and the minimal effective combination is obtained by adding two copies of each. Thus we have $$\alpha_Z = [S/Z] + 2 \cdot [S/\langle rs \rangle] + 2 \cdot [S/\langle s \rangle] $$ and its associated vector is $$\begin{pmatrix}
20 & 4 & 4 & 4 & 0 & 0 & 0 & 0 
\end{pmatrix} \ .$$

Once we have the associated vector for every $\alpha_{[Q]_\F}$, the expression of the ideals in terms of generators directly follows from Theorem \ref{thm prime ideals p-localized}.

\bibliography{biblio}
\end{document}